\documentclass[11 pt]{article}
\usepackage{pstricks,pst-plot,pst-node,pst-tree,float,multido,pstricks-add}
\usepackage{amsmath,amsfonts,graphicx,amssymb,array}

\newtheorem{defn}{Definition}
\newtheorem{thm}{Theorem}
\newtheorem{lem}{Lemma}

\newtheorem{rem}{Remark}
\newtheorem{prop}{Proposition}

\newtheorem{pf}{Proof}

\newcommand{\N}{\ensuremath{\mathbb{N}}}

\title{Avalanche polynomials}
\author{Robert Cori\thanks{LABRI, Domaine Universitaire,
 351  cours de la Lib\'eration,
 33405 Talence Cedex}, Anne Micheli and Dominique Rossin\thanks{LIAFA, Universit\'e Paris Diderot - Paris 7, Case 7014, 75205 Paris Cedex 13}}

\begin{document}
\maketitle
\begin{abstract}
The avalanche polynomial on a graph, introduced in \cite{CDR04},  capture the distribution of avalanches in the abelian sandpile model. Studied on trees, this polynomial could be defined by simply considering the size of the subtrees of the original tree.
In this article, we study some properties of this polynomial on plane trees. In \cite{CDR04}, they show that two different trees could have the same avalanche polynomial. We show here that the problem of finding a tree with a prescribed polynomial is NP-complete. In a second part, we study the average and the variance of the avalanche distribution on trees and give a closed formula.
\end{abstract}


\section{Introduction}

Self-organized criticality is a concept introduced by Bak, Tang and Wiesenfeld \cite{BTW87} to describe the behavior of natural systems like earthquakes \cite{CBO91,SS89}, forest fires. A simple model that verifies this paradigm is the Abelian Sandpile Model on the 2-dimensional lattice \cite{DRSV95,Cre91,DM94}. This model is based on a cellular automaton where each cell has a number of sand grains on it and one cell topples whenever the number of grains is greater or equal to four. In this toppling, the cell gives one grain to each of its neighbor. Thus, some other cells may topple and the sequence of toppling is called an avalanche.

This model was also considered by combinatorists \cite{Big96,BW97,CR00,GM97} but on general graphs. In \cite{CDR04}, a polynomial was introduced to encode the distribution of avalanches that was previously studied from another point of view \cite{BTW87,DM90,DV99,GLKK04,MRS02}. In this paper, we study the inverse problem. Given the distribution $D$ of avalanches, can we find a tree whose avalanche distribution is $D$ ? Moreover we study the average and the variance of this distribution on plane trees.

\section{Avalanche polynomial on a plane tree}

\subsection{Definition and main result}

The {\it avalanche polynomial} \cite{CDR04} encodes the size of avalanches in the abelian sandpile model. Let $T=(V,E,r)$ be a tree rooted at $r \in V$ whose vertex-set is $V=\{v_1,\ldots,v_{n+1}\}$, and edge-set $E$. Without loss of generality we can take $r=v_{n+1}$. Let $v \in V$ be a vertex. A {\em subtree} $T'$ rooted at $v$ is a tree whose vertices are descendants of $v$ in $T$. The {\em maximal} subtree is the subtree containing all descendants of $v$. The size of a tree is its number of vertices.

\begin{figure}[H]
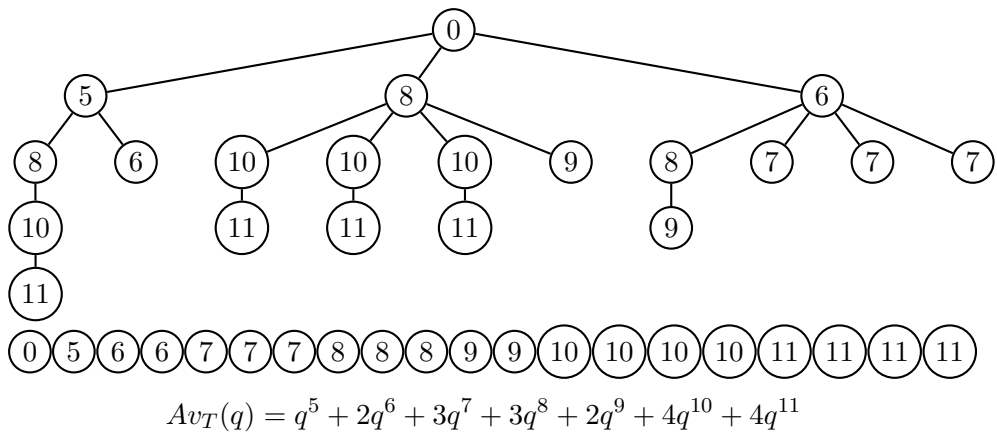

 \def\dedge{\ncline[linestyle=dotted,dotsep=1pt]}
 \psset{levelsep=25pt,labelsep=0pt,tnpos=l,radius=2pt}
\pstree{\Tcircle{0}
}{
\pstree{\Tcircle{5}
}{
\pstree{\Tcircle{8}
}{
\pstree{\Tcircle{10}
}{
\Tcircle{11}
}
}
\Tcircle{6}

}

\pstree{\Tcircle{8}
}{
\pstree{\Tcircle{10}
}{
\Tcircle{11}

}

\pstree{\Tcircle{10}
}{
\Tcircle{11}

}

\pstree{\Tcircle{10}
}{
\Tcircle{11}

}

\Tcircle{9}

}

\pstree{\Tcircle{6}
}{
\pstree{\Tcircle{8}
}{
\Tcircle{9}

}

\Tcircle{7}

\Tcircle{7}

\Tcircle{7}

}

}

\Tcircle{0} \Tcircle{5} \Tcircle{6} \Tcircle{6} \Tcircle{7} \Tcircle{7} \Tcircle{7} \Tcircle{8} \Tcircle{8} \Tcircle{8} \Tcircle{9} \Tcircle{9} \Tcircle{10} \Tcircle{10} \Tcircle{10} \Tcircle{10} \Tcircle{11} \Tcircle{11} \Tcircle{11} \Tcircle{11} 

$$Av_T(q) = q^5 + 2 q^6 + 3 q^7 + 3 q^8 + 2 q^9 + 4 q^{10} + 4 q^{11}$$
\caption{Example of an avalanche polynomial of a tree}
\label{fig:exampleLabelTree}
\end{figure}

We label the vertices $v_i$ according to the following algorithm \cite{CDR04} (see Figure \ref{fig:exampleLabelTree}):
\begin{enumerate}
\item Label the root $0$
\item The label of a child $v$ of a vertex labeled $\mu$ is $\mu + |\text{maximal subtree rooted at }v|$.
\end{enumerate}
In \cite{CDR04}, the authors show that the {\it avalanche polynomial} is also defined by $Av_T(q) = \sum_{i \geq 1} p_i q^i$, where $p_i$ is the number of vertices labeled $i$ in $T$.
In the sequel we will either speak of the polynomial or its sequence of labels.

In \cite{CDR04}, they show that the avalanche polynomial is not a tree invariant by exhibiting two different trees with the same avalanche polynomial. In the sequel, we study the problem of finding a tree for a given polyonomial.

\begin{thm}
\label{thm:thm}
Given a polynomial $P(q)$, finding whether there exists a tree $T$ whose avalanche polynomial is $P(q)$ is NP-complete.
\end{thm}

\subsection{Proof of theorem \ref{thm:thm}}

We show that the {\tt 3-PARTITION} problem can be reduced to our problem.

\begin{defn}
The {\tt 3-PARTITION}  problem $I_{n,C}$ is the following:
\begin{itemize}
\item INPUT : $C \in \N$, $(a_1,a_2,\ldots, a_{3n})$ be non-negative integers such that $\frac{C}{4} < a_i < \frac{C}{2}$.
\item OUTPUT : 
\begin{enumerate}
\item A partition of the integers into  $n$ parts of equal sum $C$.
\item NO if there is no such partition.
\end{enumerate}
\end{itemize}
\end{defn}

This is a well-known example of NP-complete problem (\cite{GJ75}).
Our proof has several steps:
\begin{enumerate}
  \item We associate a polynomial $P$ to each instance of {\tt 3-PARTITION} problem (see Definition \ref{defn:part}).
  \item We show that $P$ is an avalanche polynomial by exhibiting a tree $T_P$ with $Av_{T_P}(q) = P$. Furthermore, the structure of $T_P$ gives a solution to {\tt 3-PARTITION}. (see Lemma \ref{lem:avPol})
  \item Then, we prove that $T_P$ is the unique tree whose avalanche polynomial is $P$ (see Lemma \ref{lem:unique1}, \ref{lem:unique2}, \ref{lem:unique3}, \ref{lem:unique4}).
\end{enumerate}

\begin{defn}
\label{defn:part}
Let $I_{n,C}$ be an instance of the {\tt 3-PARTITION} problem on $a_1,\ldots,a_{3n}$.  
$P(q)$ is the polynomial associated to $I_{n,C}$:
$$P(q) = n q^{C+1} + \sum_{i=1}^{3n} q^{C+1+a_i} + \sum_{i=1}^{3n} (a_i-1) q^{C+a_i+2}$$
\end{defn}

\begin{lem}
The polynomial $P(q)$ associated to an instance $I_{n,C}$ of the {\tt 3-PARTITION} problem is an avalanche polynomial. Moreover, the tree given in Figure \ref{fig:arbreSol1} yields a solution to the {\tt 3-PARTITION} problem.
\label{lem:avPol}
\end{lem}

\begin{pf}

Let $\Pi_k=(\alpha(k),\beta(k),\gamma(k))$ be a partition of $\{1,\ldots,3n\}$ in $n$ parts 
 such that 
$$a_{\alpha(k)}+a_{\beta(k)}+a_{\gamma(k)} = C$$
Note that $P(q)$ is the avalanche polynomial of the tree $T_P$ in Figure \ref{fig:arbreSol1}.

\begin{figure}[ht]
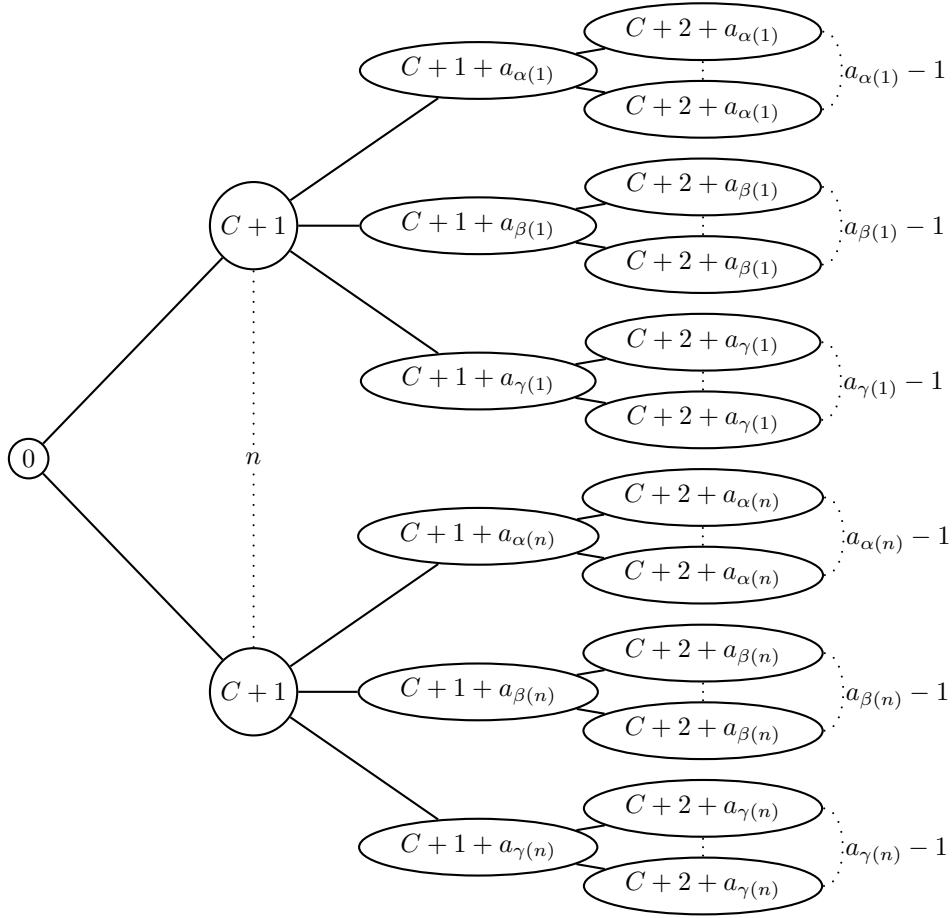

{\small
 \def\dedge{\ncline[linestyle=dotted,dotsep=1pt]}
 \psset{levelsep=85pt,treesep=.25cm,labelsep=0pt,tnpos=l,radius=0pt}
\pstree[treemode=R]{\Tcircle{0}}{
  \pstree{\Tcircle[name=n1]{$C+1$}}{
    \pstree{\Toval[name=f1]{$C+1 + a_{\alpha(1)}$}}{
      \Toval[name=pf1]{$C+2 + a_{\alpha(1)}$}
      \Toval[name=pf2]{$C+2 + a_{\alpha(1)}$}
    }
    \pstree{\Toval[name=f2]{$C+1 + a_{\beta(1)}$}}{
      \Toval[name=pf1k]{$C+2+a_{\beta(1)}$}
      \Toval[name=pf2k]{$C+2+a_{\beta(1)}$}
    }
    \pstree{\Toval[name=f3]{$C+1 + a_{\gamma(1)}$}}{
      \Toval[name=pf1h]{$C+2+a_{\gamma(1)}$}
      \Toval[name=pf2h]{$C+2+a_{\gamma(1)}$}
    }
  }
  \pstree{\Tcircle[name=n2]{$C+1$}}{
    \pstree{\Toval[name=ff1]{$C+1 + a_{\alpha(n)}$}}{
      \Toval[name=pff1]{$C+2 + a_{\alpha(n)}$}
      \Toval[name=pff2]{$C+2 + a_{\alpha(n)}$}
    }
    \pstree{\Toval[name=ff2]{$C+1 + a_{\beta(n)}$}}{
      \Toval[name=pff1k]{$C+2+a_{\beta(n)}$}
      \Toval[name=pff2k]{$C+2+a_{\beta(n)}$}
    }
    \pstree{\Toval[name=ff3]{$C+1 + a_{\gamma(n)}$}}{
      \Toval[name=pff1h]{$C+2+a_{\gamma(n)}$}
      \Toval[name=pff2h]{$C+2+a_{\gamma(n)}$}
    }
  }
}
\nccurve[linestyle=dotted]{pf1}{pf2}\trput{$a_{\alpha(1)}-1$}
\nccurve[linestyle=dotted]{pf1k}{pf2k}\trput{$a_{\beta(1)}-1$}
\nccurve[linestyle=dotted]{pf1h}{pf2h}\trput{$a_{\gamma(1)}-1$}
\ncline[linestyle=dotted]{n1}{n2}
\ncput*{$n$}
\ncline[linestyle=dotted]{pf1}{pf2}
\ncline[linestyle=dotted]{pf1k}{pf2k}
\ncline[linestyle=dotted]{pf1h}{pf2h}
\nccurve[linestyle=dotted]{pff1}{pff2}\trput{$a_{\alpha(n)}-1$}
\nccurve[linestyle=dotted]{pff1k}{pff2k}\trput{$a_{\beta(n)}-1$}
\nccurve[linestyle=dotted]{pff1h}{pff2h}\trput{$a_{\gamma(n)}-1$}
\ncline[linestyle=dotted]{pff1}{pff2}
\ncline[linestyle=dotted]{pff1k}{pff2k}
\ncline[linestyle=dotted]{pff1h}{pff2h}
}
\caption{A solution to the {\tt 3-PARTITION} problem}
\label{fig:arbreSol1}
\end{figure}

The nodes at height $1$ in $T_P$ have $3$ children, whose labels are a part of $\Pi_k$. Hence, there exists at least a solution to the {\tt Avalanche Polynomial} Problem which yields a solution to the {\tt 3-PARTITION} problem. 
\end{pf}

Nevertheless, many solutions (different non-isomorphic trees) to the {\tt Ava\-lan\-che Polynomial} Problem might exist. 

To prove  Theorem \ref{thm:thm}, we ensure the unicity  thanks to the following remark: each value $a_i$ (and $C$) can be multiplied by a polynomial integer factor $\lambda = \lambda(n)$. The avalanche polynomial becomes:
$$P_{\lambda}(q) = n q^{\lambda C+1} + \sum_{i=1}^{3n} q^{\lambda C+1+\lambda a_i} + \sum_{i=1}^{3n} (\lambda a_i-1) q^{\lambda C+\lambda a_i+2}$$
Note that $\lambda$ is polynomial ensures that our transformation is a polynomial-time reduction.
With the following substitutions ($a_i \rightarrow \lambda a_i$ and $C \rightarrow \lambda C$), the tree of Figure \ref{fig:arbreSol1} is a solution of $P_{\lambda}$.
We are going to prove that, taking $\lambda$ large enough, every tree such that $Av_T(q) = P_{\lambda}(q)$ yields a solution to the {\tt 3-PARTITION} problem.

\begin{lem}\label{lem:unique1}
The root has $n$ children, each one labeled {$\lambda C+1$}.
\end{lem}

\begin{pf}
By definition, the smallest value is the label of a child of the root. Hence the $n$ vertices labeled  $\lambda C+1$ are children of the root.
But $n (\lambda C + 1) = (\sum  \lambda a_i) + n$, the number of vertices. Hence all other vertices are children of one of the vertices labeled  $\lambda C + 1$.
\end{pf}

Let $I$ be the set of all vertices labeled $\lambda C + 1 + \lambda
a_i$ and $L$ the set of those labeled $\lambda C + \lambda a_i+2$,
$\lambda \geq 2$. Note that if $\lambda > 1$ then $|I| = 3n$, $|L| = n(\lambda C -3)$ and $I \bigcap L = \varnothing$.

\begin{lem}
\label{lem:feuille1}\label{lem:unique2}
If $\lambda > 2$ then the leaves of every solution tree are in $L$.
\end{lem}

\begin{pf}
Let $v$ be a leaf of a solution tree labeled $l_v$. Then its parent $p$ is labeled $l_v-1$. If $l_v = \lambda C + 1 + \lambda a_v$ there does not exist any vertex labeled $\lambda C + \lambda a_v$. Otherwise $\lambda C +\lambda a_v  = \lambda C + \lambda a_p + 1 + \delta$ ($\delta \in \{0,1\}$). Hence, $\lambda (a_v-a_p) = 1 + \delta$ which contradicts $\lambda > 2$. $p$ cannot be labeled $\lambda C+1$ since $\lambda a_v \neq 1$.

\end{pf}

\begin{lem}
\label{lem:feuille2}\label{lem:unique3}
If $\lambda > 3n$ then all vertices of $L$ are leaves.
\end{lem}

\begin{pf}
Suppose that there exists an internal vertex in $L$. Take $v$ as the one with the greatest label $l_v = \lambda C + \lambda a_v + 2$.
 Let $\nu$ be a child of $v$. If $\nu$ is a leaf, then its label is $l_v+1=\lambda C+\lambda a_{v}+3$. But there is no such vertex since $\lambda > 3$. 
Hence $\nu$ is an internal node. But $v$ is the greatest internal node of $L$ and the label of $\nu$ is greater than $l_v$, hence $\nu \in I$.

Then $v$ is either a child of the node labeled $\lambda C +1$ or not. The two cases are shown in Figure \ref{fig:shape}.

\begin{figure}[ht]
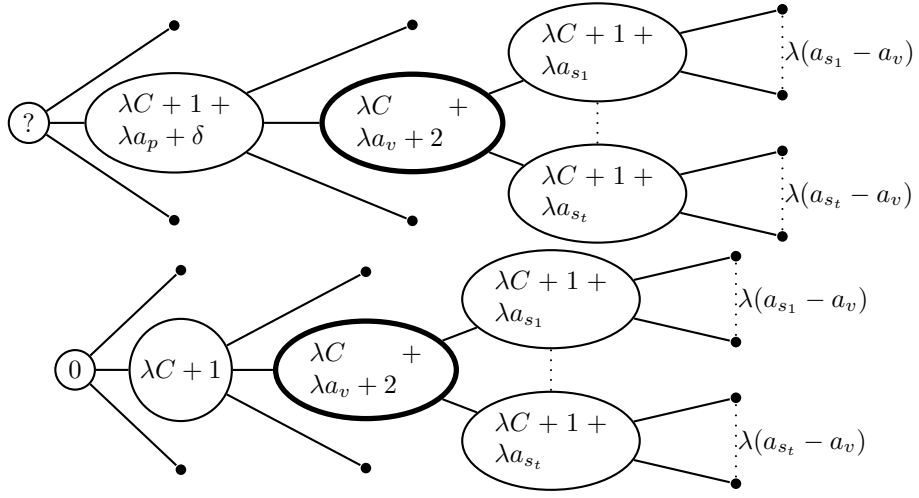

\begin{center}
{\small
  \def\dedge{\ncline[linestyle=dotted,dotsep=1pt]}
  \psset{levelsep=55pt,treesep=.55cm,labelsep=0pt,tnpos=l,radius=0pt}

\pstree[treemode=R,levelsep=55pt]{\Tcircle{?}}{
  \Tc*{2pt}
  \pstree[treemode=R,levelsep=90pt]{\Toval{\parbox[t]{1.5cm}{$\lambda C + 1 + \lambda a_p + \delta$}}}{
    \Tc*{2pt}
    \pstree[levelsep=70pt]{\Toval[name=labelvv,linewidth=2pt]{\parbox[t]{1.5cm}{$\lambda C + \lambda a_v + 2$}}}{
      \pstree[levelsep=70pt,treesep=1cm]{\Toval[name=pfff1]{\parbox[t]{1.5cm}{$\lambda C + 1+\lambda a_{s_1}$}}}{
	\Tc*[name=fann1]{2pt}
	\Tc*[name=fann1p]{2pt}
      }
      \pstree[levelsep=70pt,treesep=1cm]{\Toval[name=pfff2]{\parbox[t]{1.5cm}{$\lambda C + 1+\lambda a_{s_t}$}}}{
	\Tc*[name=fann2]{2pt}
	\Tc*[name=fann2p]{2pt}

      }
    }
    \Tc*{2pt}

  }
  \Tc*{2pt}
}
  
\ncline[linestyle=dotted]{fann1}{fann1p}\trput{$\lambda(a_{s_1}-a_v)$}
\ncline[linestyle=dotted]{fann2}{fann2p}\trput{$\lambda(a_{s_t}-a_v)$}
\ncline[linestyle=dotted]{pfff1}{pfff2}

  \pstree[levelsep=40pt,treemode=R]{\Tcircle{0}}{
    \Tc*{2pt}
    \pstree[levelsep=70pt]{\Tcircle{$\lambda C + 1$}}{
      \Tc*{2pt}
      \pstree{\Toval[name=labelv,linewidth=2pt]{\parbox{1.5cm}{$\lambda C  + \lambda a_v + 2$}}}{
	\pstree[levelsep=70pt,treesep=1cm]{\Toval[name=pff1]{\parbox{1.5cm}{$\lambda C + 1+\lambda a_{s_1}$}}}{
	\Tc*[name=fan1]{2pt}
	\Tc*[name=fan1p]{2pt}
	}
	\pstree[levelsep=70pt,treesep=1cm]{\Toval[name=pff2]{\parbox{1.5cm}{$\lambda C + 1+\lambda a_{s_t}$}}}{
	\Tc*[name=fan2]{2pt}
	\Tc*[name=fan2p]{2pt}
	}
      }
      \Tc*{2pt}
    }
    \Tc*{2pt}
  }
  
\ncline[linestyle=dotted]{fan1}{fan1p}\trput{$\lambda(a_{s_1}-a_v)$}
\ncline[linestyle=dotted]{fan2}{fan2p}\trput{$\lambda(a_{s_t}-a_v)$}
\ncline[linestyle=dotted]{pff1}{pff2}
}
\caption{Two different shapes of trees that appear in Lemma \ref{lem:feuille2}.}
\label{fig:shape}
\end{center}
\end{figure}

Suppose that the parent $p$ of $v$ is labeled $\lambda C + 1 + \lambda a_p + \delta, \delta = 0,1$.
The size $\sigma$ of the subtree rooted at $v$ is given:
\begin{enumerate}
\item either by the difference between its label and the label of its parent,

$$\sigma  =  \lambda C + \lambda a_v +2 - (\lambda C + 1 + \lambda a_p + \delta)\\
  =  \lambda(a_v-a_p) + (1-\delta),$$

\item or by the sum of the size of the $t$ child-trees of $v$ plus $1$ (to take $v$ into account):
$$\sigma =  1 + \sum_{i = 1}^{t} \left( (\lambda C + \lambda a_{s_i} + 1) - (\lambda C+\lambda a_v + 2) \right) 
= 1 - t + \sum_{i=1}^{t} \lambda (a_{s_i} - a_v)$$
\end{enumerate}

Hence 
$\lambda(a_v-a_p) + (1 - \delta) = 1 - t + \sum_{i=1}^{t} \lambda (a_{s_i} - a_v)$, and\\
$\lambda \left(\sum_{i=1}^{t} a_{s_i} - (t+1)a_v + a_p \right) = t - \delta $.\\
But $t \leq |I| = 3n$ (the number of values $a_i$). Hence $t-\delta < 3n$. Taking $\lambda > 3n$ concludes whenever $p$ is labeled $\lambda C + \lambda a_p + 1 + \delta, \delta = 0,1$.

Note that when $p$ is labeled $\lambda C + 1$, taking $a_p = 0$ and $\delta = 0$ in the preceding proof concludes.

\end{pf}

\begin{lem}\label{lem:unique4}
If $\lambda > 3n$, then the nodes of $I$ are children of one of the vertices labeled $\lambda C + 1$.
\end{lem}
\begin{pf}
Let $v \in I$ be a vertex labeled $\lambda C +\lambda a_v + 1$. By
Lemmas \ref{lem:feuille1} and \ref{lem:feuille2}, the $\lambda a_v -1$
vertices labeled $\lambda C +\lambda a_v + 2$ are leaves, hence
children of $v$. Note that if $k$ vertices $v_1,\ldots,v_k$ have the
same label than $v$ then there are $k (\lambda a_v -1)$ vertices
labeled $\lambda C +\lambda a_v + 2$. These vertices are children of
the $v_i, 1 \leq i \leq k$. The size of each maximal subtree rooted at $v_i$ is less or
equal than $\lambda a_v$ since the size of the subtree is the label of
$v$ minus the label of its parent $p$, and the label of $p$ is at
least $\lambda C+1$.  Hence each $v_i$ has exactly $\lambda a_v -1$
children labeled $\lambda C + \lambda a_v+2$.

Suppose that the parent $p$ of $v$ is labeled $\lambda C + 1 + \lambda a_p $ then the size of the subtree rooted at $v$ is $\lambda (a_v-a_p) < \lambda a_v$.

Thus the parent of $v$ is labeled $\lambda C + 1$.
\end{pf}


\subsection{Avalanche polynomial on plane tree of height at most $2$.}

We showed in the last section that finding a tree whose avalanche polynomial is given is NP-complete. In this part, we show that considering only trees of height at most $2$, then the problem becomes poynomial. Moreover we give a linear time and space algorithm to solve it.

\begin{figure}[ht]
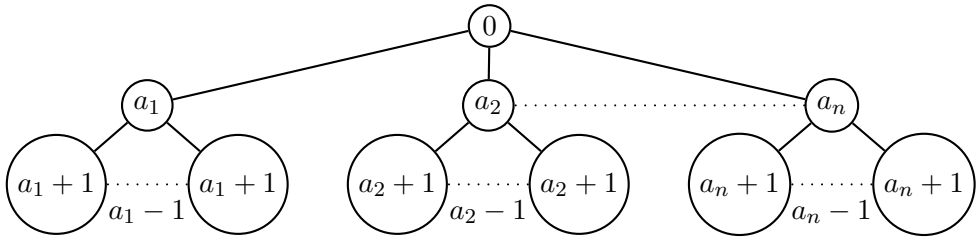

\begin{center}
\psset{levelsep=30pt}
\pstree{
  \Tcircle{0}}{
    \pstree[treesep=30pt]{
      \Tcircle{$a_1$}}{
	\Tcircle[name=a11]{$a_1+1$}
	\Tcircle[name=a12]{$a_1+1$}
      }
    
    \pstree[treesep=30pt]{
      \Tcircle[name=a2]{$a_2$}}{
	\Tcircle[name=a21]{$a_2+1$}
	\Tcircle[name=a22]{$a_2+1$}
      }
    
    \pstree[treesep=30pt]{
      \Tcircle[name=an]{$a_n$}}{
	\Tcircle[name=an1]{$a_n+1$}
	\Tcircle[name=an2]{$a_n+1$}
      }
    }
\ncline[linestyle=dotted]{a2}{an}
\ncline[linestyle=dotted]{a11}{a12}\tbput{$a_1-1$}
\ncline[linestyle=dotted]{a21}{a22}\tbput{$a_2-1$}
\ncline[linestyle=dotted]{an1}{an2}\tbput{$a_n-1$}
\caption{Labeled tree of height $2$}
\label{fig:arbre2}
\end{center}
\end{figure}

A tree of height $2$ and its labeling is given in Figure \ref{fig:arbre2} and its associated polynomial is $P(q) = \sum_{i=1}^{n} \left( q^{a_i} + (a_i-1) q^{a_i+1} \right)$. Note that if one of the $a_i$ is equal to $1$ then it is a leaf rooted at the vertex labeled $0$.

The following algorithm takes as input a polynomial $P(q) = \sum_{i=1}^{n} a_i q^i$ and output either:
\begin{itemize}
\item A tree $T$ of height at most $2$ whose avalanche polynomial is $P$.
\item NO if no such tree exists.
\end{itemize}

\begin{enumerate}
\item Create the root with label $0$. 
\item \label{step2}Find a child of the root by taking the first non-nul coefficient $a_{j}$. If $j = 1$ we have found $a_j$ leaves of the root. If $j > 1$ we have $a_j$  nodes. Each of these nodes must have $j-1$ children labeled $j+1$. Hence $a_{j+1}$ must be greater than $a_j (j-1)$. 
If not, output NO and exit. 
\item Subtract $a_j q^j + a_j (j-1) q^{j+1}$ from the polynomial;  we take this new polynomial as $P$
\item If $P = 0$ output the tree and exit.
\item Goto \ref{step2}.
\end{enumerate}

\section{Avalanche polynomial on plane trees}

In this section, we study the avalanche polynomial on plane trees. The aim of this section is to find the average distribution of avalanches on a random tree. We first give the closed expression for the number of avalanches of a given size in plane trees. Unfortunately, this formula does not help to retrieve informations on the distribution, thus  we give a recursive formula for the distribution and deduce the mean and the variance of this distribution.

\subsection{Definition}

Let us note $C_k = \frac{1}{k+1}\binom{2k}{k}$ the $k$-th Catalan number. Let ${\mathcal T}_n$ be the set of rooted plane trees with $n$ edges. $|{\mathcal T}_n| = C_n$. 
We extend the definition of avalanche polynomial to the set of plane trees (see Figure \ref{fig:tableau}):
\begin{equation}
A(t,q)  =  \sum_{i \geq 0} \sum_{T \in {\mathcal T_i}} Av_{T}(q) t^i  =  \sum_{p\geq 0} A_p(q) t^p
\label{eq:atq}
\end{equation}

\begin{figure}[ht]
\begin{center}
\begin{tabular}{|c|c|c|}
\hline
Size & Trees & Polynomial\\
\hline
1 & {\tiny \psset{levelsep=0pt,treesep=.25cm,labelsep=0pt,tnpos=l,radius=0pt}
\pstree{\Tcircle{0}}{}
} & $0$\\
\hline
2 & 
{\tiny \psset{levelsep=20pt,treesep=.25cm,labelsep=0pt,tnpos=l,radius=0pt}
\pstree{\Tcircle{0}}{
  \Tcircle{1}}
}
& $qt$\\
\hline
3 & 
{\tiny \psset{levelsep=20pt,treesep=.25cm,labelsep=0pt,tnpos=l,radius=0pt}
\pstree{\Tcircle{0}}{
  \Tcircle{1}
  \Tcircle{1}}

}
{\tiny \psset{levelsep=20pt,treesep=.25cm,labelsep=0pt,tnpos=l,radius=0pt}
\pstree{\Tcircle{0}}{
  \pstree{\Tcircle{2}}{
  \Tcircle{3}}
}}
& $(2q+q^2+q^3)t^2$\\
\hline
4 &
{\tiny \psset{levelsep=20pt,treesep=.25cm,labelsep=0pt,tnpos=l,radius=0pt}
\pstree{\Tcircle{0}}{
  \Tcircle{1}
  \Tcircle{1}
  \Tcircle{1}}

}
{\tiny \psset{levelsep=20pt,treesep=.25cm,labelsep=0pt,tnpos=l,radius=0pt}
\pstree{\Tcircle{0}}{
  \Tcircle{1}
  \pstree{\Tcircle{2}}{
  \Tcircle{3}}
}}
{\tiny \psset{levelsep=20pt,treesep=.25cm,labelsep=0pt,tnpos=l,radius=0pt}
\pstree{\Tcircle{0}}{
  \pstree{\Tcircle{2}}{
  \Tcircle{3}}
  \Tcircle{1}
}}
{\tiny \psset{levelsep=20pt,treesep=.25cm,labelsep=0pt,tnpos=l,radius=0pt}
\pstree{\Tcircle{0}}{
  \pstree{\Tcircle{3}}{
  \Tcircle{4}
  \Tcircle{4}}

}}
{\tiny \psset{levelsep=20pt,treesep=.25cm,labelsep=0pt,tnpos=l,radius=0pt}
\pstree{\Tcircle{0}}{
  \pstree{\Tcircle{3}}{
  \pstree{\Tcircle{5}}{
  \Tcircle{6}}}

}}
& $\begin{array}{l} 
    (5q+2q^2+4q^3 \\
    +2q^4+q^5+q^6)t^3
  \end{array}$\\
\hline

\end{tabular}
\caption{First terms of $A(t,q)$}
\label{fig:tableau}
\end{center}
\end{figure}

\begin{rem}
\label{rk:t1}
 Note that :
\begin{equation*}
A(t,1) = \sum_{p\geq 0} A_p(1) t^p = \sum_{p\geq 0} p C_p t^p
\end{equation*}
\end{rem}

In fact, the coefficient of $t^p$ in $A(t,1)$ is the number of rooted plane trees of size $p+1$ times the number of vertices (excepting the roots). Hence it is $p C_p$.

\subsection{Closed formula}

In this section, we exhibit a decomposition of the rooted plane trees which yields a closed formula for the avalanche polynomial.

\begin{thm}
The coefficients of the polynomials $A_n(q)$ are given by:
\begin{equation}
[q^v] A_n(q) = \sum_{k = 1}^{{\left\lfloor{\frac{-1+\sqrt{1+8v}}{2}}\right\rfloor}} \sum_{p \in \Pi_{n,v,k}} C_{p_1-1} \left( \prod_{i=2}^{k} C_{p_i-p_{i-1}} \right) C_{n-p_k+1} 
\label{eq:thm}
\end{equation}

where 
$$\Pi_{n,v,k} = \{ (p_1,\ldots,p_k); p_i < p_{i+1},  p_k \leq n, \sum_{i=1}^{k} p_i = v \} $$
\end{thm}

\begin{pf}

The coefficient $[q^v] A_n(q)$ represents the number of vertices labeled $v$ in every plane trees of size $n+1$.

Let $x$ be such a vertex labeled $v$ in a tree $T$ and $y_0=x,y_1,\ldots,y_k = x_{n+1}$ be the path joining the root $x_{n+1}$ labeled $0$ to the vertex $x$ in $T$. 
Let $lab$ be the function which maps any vertex $y$ of $T$ onto $lab(y)$, the label of $y$ in $T$.
The proof is based on the study of rooted plane trees which have a given path.

\begin{itemize}
\item Note that $k \leq \left\lfloor \frac{-1+\sqrt{1+8v}}{2} \right\rfloor$ because the labels of the vertices of the path from the root to $x$ take minimal values when the maximal subtree of $T$  rooted at $y_{k-1}$ is reduced to the path $y_0,\ldots,y_{k-1}$.

Then the labels of the vertices on this path are $k$, $k+(k-1)$,$k+(k-1)+(k-2)$, ..., $\frac{k(k+1)}{2}$ and if $k > \left\lfloor \frac{-1+\sqrt{1+8v}}{2} \right\rfloor$ then the last label is stricly greater than $v$.
\item The sequence of labels of the  vertices $y_l$ is stricly decreasing by definition.
\item The sequence of the differences $(lab(y_l)-lab(y_{l+1}))$ is
stricly increasing. This relies on the fact that each difference is
the size of the maximal subtree rooted at $y_{l+1}$ and that the
subtrees are nested. Moreover, $lab(y_{k-1})-lab(y_{k}) \leq n$ as this is the
size of the maximal subtree rooted at $y_{k-1}$. Finally,
$\sum_{i=1}^{k} lab(y_{i-1})-lab(y_{i}) = lab(y_0) - lab(y_k)= v$.

\end{itemize}

Let $k$ be an integer such that $k \leq \left\lfloor \frac{-1+\sqrt{1+8v}}{2} \right\rfloor$ and let $(p_1,\ldots,p_k)$ be a stricly increasing sequence such that $p_k \leq n$ and $\sum_{i=1}^{k} p_i = v$. Let $P$ be a path of vertices $\{y_0,\ldots,y_k\}$ labeled $lab(y_0)=v, lab(y_1) = lab(y_0)-p_1,  lab(y_2) = lab(y_1)-p_2, \ldots, lab(y_k) =  lab(y_{k-1}) - p_k = 0$ (see Figure \ref{fig:P}). Denote by $F_{i,1}$ (resp. $F_{i,2}$) the forest hanging on $y_i$ at the left (resp. right) of the edge $(y_i,y_{i-1})$.

\begin{lem}
The number of plane trees of size $n+1$ having $P$ as subpath is $C_{p_1-1}\left( \prod_{i=2}^{k} C_{p_i-p_{i-1}} \right) C_{n-p_k+1}$.
\end{lem}
\begin{pf}

Note first that the size of the maximal subtree hanging on $y_i$ is $y_i-y_{i+1} = p_{i+1}$. Then the number of vertices of the forest $F_{i,1} \cup F_{i,2}$ is $p_{i+1}-p_i-1$ for $i \in \{ 1,\ldots,k-1 \}$. The number of such pair of forests $\{F_{i,1},F_{i,2}\}$ is $C_{p_{i+1}-p_i}$.

The size of the maximal subtree rooted at $y_0$ is $y_0-y_1=p_1$. There are $C_{p_1-1}$ such trees.

The remaining nodes are in the  forests $F_{k,1}$ and $F_{k,2}$ which union has $n-p_k$ vertices. So that the number of such forests is $C_{n-p_k+1}$.

\begin{figure}[H]
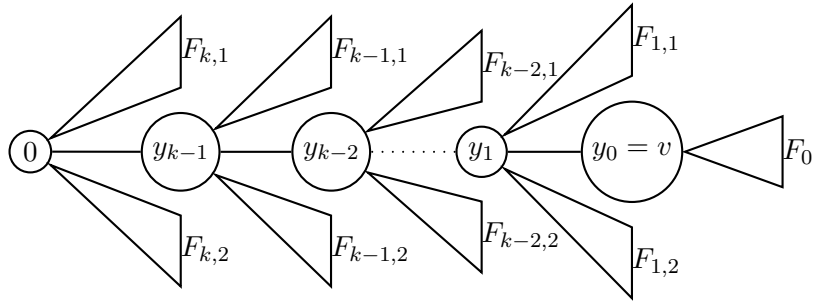

\begin{center}
\pstree[treemode=R,treesep=.3cm]{\Tcircle{$0$}}{ \Tfan[name=x01]
  \pstree[linestyle=solid]{\Tcircle{$y_{k-1}$}}{
  \Tfan[name=x11,linestyle=solid]
  \pstree[linestyle=dotted]{\Tcircle{$y_{k-2}$}}{
  \Tfan[name=x21,linestyle=solid]
  \pstree[linestyle=solid]{\Tcircle[linestyle=solid]{$y_{1}$}}{ \Tfan[name=xk1]
  \pstree{\Tcircle{$y_{0}=v$}}{ \Tfan[name=fin] } \Tfan[name=xk2]

    }
    \Tfan[name=x22,linestyle=solid]          
    }
    \Tfan[name=x12,linestyle=solid]
  }
  \Tfan[name=x02]
}
\rput[l](x01){$F_{k,1}$}
\rput[l](x02){$F_{k,2}$}
\rput[l](x11){$F_{k-1,1}$}
\rput[l](x12){$F_{k-1,2}$}
\rput[l](x21){$F_{k-2,1}$}
\rput[l](x22){$F_{k-2,2}$}
\rput[l](xk1){$F_{1,1}$}
\rput[l](xk2){$F_{1,2}$}
\rput[l](fin){$F_0$}
\caption{A tree with its subpath $P$}
\label{fig:P}
\end{center}
\end{figure}

\end{pf}

\end{pf}

\begin{rem}
The coefficient of $t^nq$ in $A(t,q)$ is $C_{n}$ and the coefficient of $t^nq^2$ is $C_{n-1}$.
\label{rk:Terme1}
\end{rem}

\begin{pf}
There are two different proofs of these results, the first one consists in taking $v=1$ or $v=2$ in Equation (\ref{eq:thm}).

The second one relies on the following property:

\begin{lem}
Let $T'$ be a tree. Let $S_{T'}$ be the set of trees $T$ of size $n+1$ where $T'$ is a subtree of $T$  rooted at the root of $T$. Then $|S_{T'}| = C_{n-|T'|+2}$.
\end{lem}

\begin{pf}
The bijection between the trees of $S_{T'}$ and plane trees of size $n-|T'|+1$  is given in Figure \ref{fig:bijectionArbre}. 

\begin{figure}
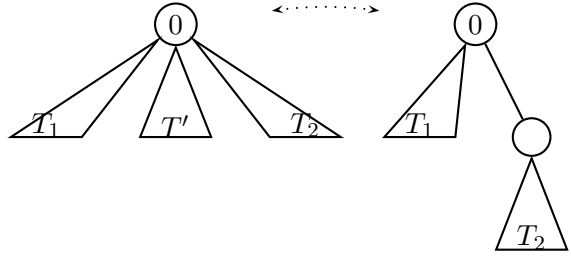

\begin{center}

\begin{tabular}{lr}

\pstree[levelsep=1.5cm]{\Tcircle[name=r1]{0}}{
\Tfan[name=T1]
\Tfan[name=Tp]
\Tfan[name=T2]
}

\rput[b](T1){$T_1$}
\rput[b](Tp){$T'$}
\rput[b](T2){$T_2$}

&
\pstree[levelsep=1.5cm]{\Tcircle[name=r2]{0}}{
\Tfan[name=TT1]
\pstree{\TC}{
\Tfan[name=TT2]
}
}
\rput[b](TT1){$T_1$}
\rput[b](TT2){$T_2$}
\ncarc[nodesep=1cm,linestyle=dotted]{<->}{r1}{r2}
\end{tabular}
\caption{Bijection between $S_{T'}$ and plane trees}
\label{fig:bijectionArbre}
\end{center}
\end{figure}
\end{pf}

For $v=1$ take $T'$ as a single edge. For $v=2$ take $T'$ as a path of length $2$.

\end{pf}

\subsection{Functional equation}

\begin{lem}
The polynomial $A_{p}(q)$, $p\geq 0$ - the distribution of avalanches on plane trees of size $p+1$- is defined by the following recurrence:
\begin{equation}
A_{p+1}(q) = \sum_{k=0}^p C_k C_{p-k} q^{k+1} + C_{p-k} q^{k+1} A_k(q) + C_k A_{p-k}(q)
\label{eq:ap1}
\end{equation}
This yields a functionnal equation for $A(t,q)$:
\begin{eqnarray}
A(t,q)&=& \frac{(1-\sqrt{1-4t})(1-\sqrt{1-4qt})}{4t} \nonumber\\
& &+q\frac{1-\sqrt{1-4t}}{2}A(qt,q)+\frac{1-\sqrt{1-4t}}{2}A(t,q)
\label{eq:atq2}
\end{eqnarray}

\end{lem}

\begin{pf}

The proof relies on the following decomposition  of rooted plane trees (see Figure \ref{fig:decomposition}). A tree $T_{p+1}$ of size $p+2$ can be decomposed into a tree $T_k$ of size $k+1$ and a tree $T_{p-k}$ of size $p+1-k$ as shown in Figure \ref{fig:decomposition}. Hence, the root of the tree $T_k$ is the leftmost child of the root of $T_{p-k}$ in $T_{p+1}$. Furthermore, the root of $T_k$ is labeled $k+1$ the number of vertices of $T_k$. This vertex  contributes to the term $C_k C_{p-k} q^{k+1}$ in $A_{p+1}(q)$. 

The avalanche polynomial of $T_k$ is $A_k(q)$. In $T_{p+1}$ the label of a vertex of the subtree $T_k$ is increased by $k+1$. They contribute to the term $C_{p-k} q^{k+1} A_k(q)$ in $A_{p+1}(q)$. The last term takes into account the labels of the vertices of $T_{p-k}$ which stay unchanged since the root of $T_{p-k}$ is the root of $T_{p+1}$.

\begin{figure}[H]
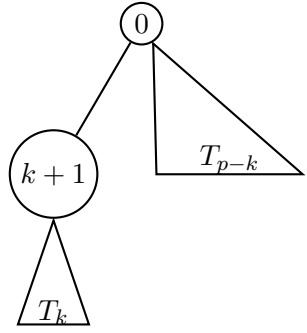

\begin{center}
\pstree{\Tcircle{$0$}}{
  \pstree{\Tcircle{$k+1$}}{
    \Tfan[name=Tk]}
  \Tfan[fansize=2,name=Tpk]}
\rput[b](Tk){$T_k$}
\rput[b](Tpk){$T_{p-k}$}
\caption{Decomposition of a plane tree}
\label{fig:decomposition}
\end{center}
\end{figure}

Substituting Equation (\ref{eq:ap1}) in Equation (\ref{eq:atq}) yields a functionnal equation for $A(t,q)$:
\begin{eqnarray*}
A(t,q) &= &\frac{(1-\sqrt{1-4t})(1-\sqrt{1-4qt})}{4t} \nonumber\\
 & &+ \sum_{p\geq 0}{\sum_{k=0}^p (C_{p-k} q^{k+1} A_k(q) t^{p+1} + C_k A_{p-k}(q) t^{p+1})} \nonumber\\
&=& \frac{(1-\sqrt{1-4t})(1-\sqrt{1-4qt})}{4t} \nonumber\\
& &+q\frac{1-\sqrt{1-4t}}{2}A(qt,q)+\frac{1-\sqrt{1-4t}}{2}A(t,q)
\end{eqnarray*}
\end{pf}



\subsection{Average of the distribution of avalanches}

\begin{prop}
The average size of an avalanche on a plane tree of size $n$ is asymptotically 
equal to $\frac{\sqrt \pi}{4} n^{3/2}$.
\end{prop}

\begin{pf}

$$A(t,q) = \sum_n \sum_m a_{m,n} q^m t^n$$

The mean size of avalanches on all plane trees of size $n$ is given by 
\begin{equation}
M_n = \frac{\sum_m{m a_{m,n}}}{\sum_m a_{m,n}}  = \frac{\left[ \frac{\partial A}{\partial q}(t,1)\right]_n}{\left[ A(t,1) \right]_n} = \frac{\left[ \frac{\partial A}{\partial q}(t,1)\right]_n}{n C_n}
\label{eq:Mn}
\end{equation}

Equation (\ref{eq:atq2}) yields:
\begin{eqnarray*}
\frac{\partial A}{\partial q}(t,1) & = & \frac{1-\sqrt{1-4t}}{2\sqrt{1-4t}}+\frac{1-\sqrt{1-4t}}{2} \left( A(t,1) + \frac{\partial A}{\partial q}(t,1) +t\frac{\partial A}{\partial t}(t,1) + \frac{\partial A}{\partial q}(t,1) \right) 
\end{eqnarray*}

Factorization of  $\frac{\partial A}{\partial q}(t,1)$ and substitution of  $A(t,1)$ and $\frac{\partial A}{\partial t}(t,1)$ with their values of Remark \ref{rk:t1} yield:

\begin{eqnarray*}
\frac{\partial A}{\partial q}(t,1)& = & \frac{(\sqrt{1-4t}-1)(2t-1)}{2(4t-1)^2}
\end{eqnarray*}

Since
\begin{equation*}
\frac{(\sqrt{1-4t}-1)(2t-1)}{2}  =   t + \sum_{i\geq 0} (C_{i+1}-2C_i)t^{i+2},
\end{equation*}

\begin{eqnarray*}
\frac{\partial A}{\partial q}(t,1) & = & t\sum_{i\geq 0} 4^i (i+1) t^i + \sum_{i,j\geq 0}(C_{i+1}-2C_i) 4^j (j+1) t^{i+j+2}\\
& = & t\sum_{i\geq 0} 4^i (i+1) t^i + \sum_{k\geq 1} 4^k t^{k+1} \sum_{i=0}^{k-1} (C_{i+1}-2C_i) \frac{k-i}{4^{i+1}}
\end{eqnarray*}

\noindent Since
$\displaystyle\sum_{i=0}^{k-1} \frac{C_{i+1}-2C_i}{4^{i+1}} =  -\frac{k C_k}{4^k},
\sum_{i=0}^{k-1} \frac{i (C_{i+1}-2C_i)}{4^{i+1}}  = \frac{(2+5k+k^2) C_k}{4^k}-2$,

\begin{eqnarray*}
\frac{\partial A}{\partial q}(t,1) & = & \sum_{i\geq 1} \left(4^{i-1} (i+2) - (-1+i+2i^2)C_{i-1}\right)t^i
\end{eqnarray*}

\noindent Substitution of $\frac{\partial A}{\partial q}(t,1)$ in Equation (\ref{eq:Mn}) yields:
\begin{equation*}
M_n  = \frac{4^{n-1}(n+2) }{n C_n}-\frac{(n+1)^2}{2n}
\end{equation*}

\noindent Since $C_n \sim \frac{4^n}{\sqrt{\pi n^3}}$,
$M_n \sim \frac{\sqrt{\pi}}{4} n^{3/2}$.

\end{pf}

\subsection{Variance}

The variance of the distribution is given by 
$\frac{\left[ \frac{\partial^2 A}{\partial q^2}(t,1) \right]_n}{n C_n} + M_n - M_n^2$ where $M_n$ denotes the mean size (computed in the last section).

Hence we only need the expression of $\left[ \frac{\partial^2 A}{\partial q^2}(t,1) \right]_n$.

\begin{equation*}
\begin{split}
\frac{\partial^2 A}{\partial q^2}(t,1) &= \frac{t(1-\sqrt{1-4t})}{(1-4t)^{3/2}} + \frac{1-\sqrt{1-4t}}{2} \left( 2t \frac{ \partial A}{\partial t}(t,1) + 2 \frac{\partial A}{\partial q} (t,1)   \right.\\
& \quad \left. + t^2 \frac{\partial^2 A}{\partial t^2}(t,1) + 2 \frac{\partial^2 A}{\partial q^2}(t,1) + 2 t \frac{\partial^2 A}{\partial t\partial q}(t,1) \right)
\end{split}
\end{equation*}

Thus we need the other derivatives $\frac{\partial A}{\partial t}(t,q)$,$\frac{\partial A}{\partial q}(t,q)$, $\frac{\partial^2 A}{\partial t^2}(t,1)$ and $\frac{\partial^2 A}{\partial t \partial q}(t,1)$.
Substituting their expression in $\frac{\partial^2 A}{\partial q^2}(t,1)$ yields:

\begin{equation*}
\begin{split}
\frac{\partial^2 A}{\partial q^2}(t,1) &=
-1/2\,{\frac { \left( -1+\sqrt {1-4\,t} \right)  \left( 10\,t+\sqrt {1
-4\,t}-1-32\,{t}^{2}+32\,{t}^{3} \right) }{ \left( -1+4\,t \right) ^{4
}}}
\end{split}
\end{equation*}


Since :
\begin{equation*}
\begin{split}
(1-4t)^{-4} &= \sum_{n\geq 0} \frac{4^n}{6} (n+1)(n+2)(n+3) t^n\\
(1-4t)^{-7/2} &= \sum_{n \geq 0} \frac{(n+1) (2n+1)(2n+3)(2n+5) 4^n C_n}{15 2^{2n}}t^n
\end{split}
\end{equation*}
We obtain:
\begin{equation*}
\begin{split}
V_n &= \frac{4}{15}n^3+\frac{73}{60}n^2+\frac{26}{15}n+\frac{8}{15}-\frac{1}{2n}-\frac{1}{2n}-\frac{1}{4n^2}\\
& \quad -\frac{16^{n-1}}{n^2(n+1)^2C_n^2}(n^4+6n^3+13n^2+12n+4)\\
& \quad +\frac{4^{n-1}}{n^2(n+1)C_n}(n^3+4n^2+5n+2)
\end{split}
\end{equation*}

The leading term is:
$$\frac{4}{15}n^3-\frac{16^{n-1}}{n^2(n+1)^2C_n^2}n^4 \sim \left( \frac{4}{15} - \frac{\pi}{16} \right) n^3$$

In fact the distribution show strange peaks and after a renormalization (scaling in $x$ by a factor $1/n$ for the distribution on trees of size $n$ and rescaling the ordinate between $0$ and $1$) it gives the curves in Figure \ref{fig:des1} and \ref{fig:des2}. In the following figures, we draw a curve for $A_{10}(q)$, $A_{20}(q)$,... Moreover for such a polynomial $A_n(q) = \sum_{i \geq 0} p_i q^i$ we put points at coordinates $(i/n, p_i/p_1)$ since $p_1$ is the maximum. Let's recall from Remark \ref{rk:Terme1} that $p_1 = C_n$. The grey scale goes from white (for the first polynomial) to black (for the last one).

\begin{figure}[H]
\begin{center}
\includegraphics[width=.45\textwidth]{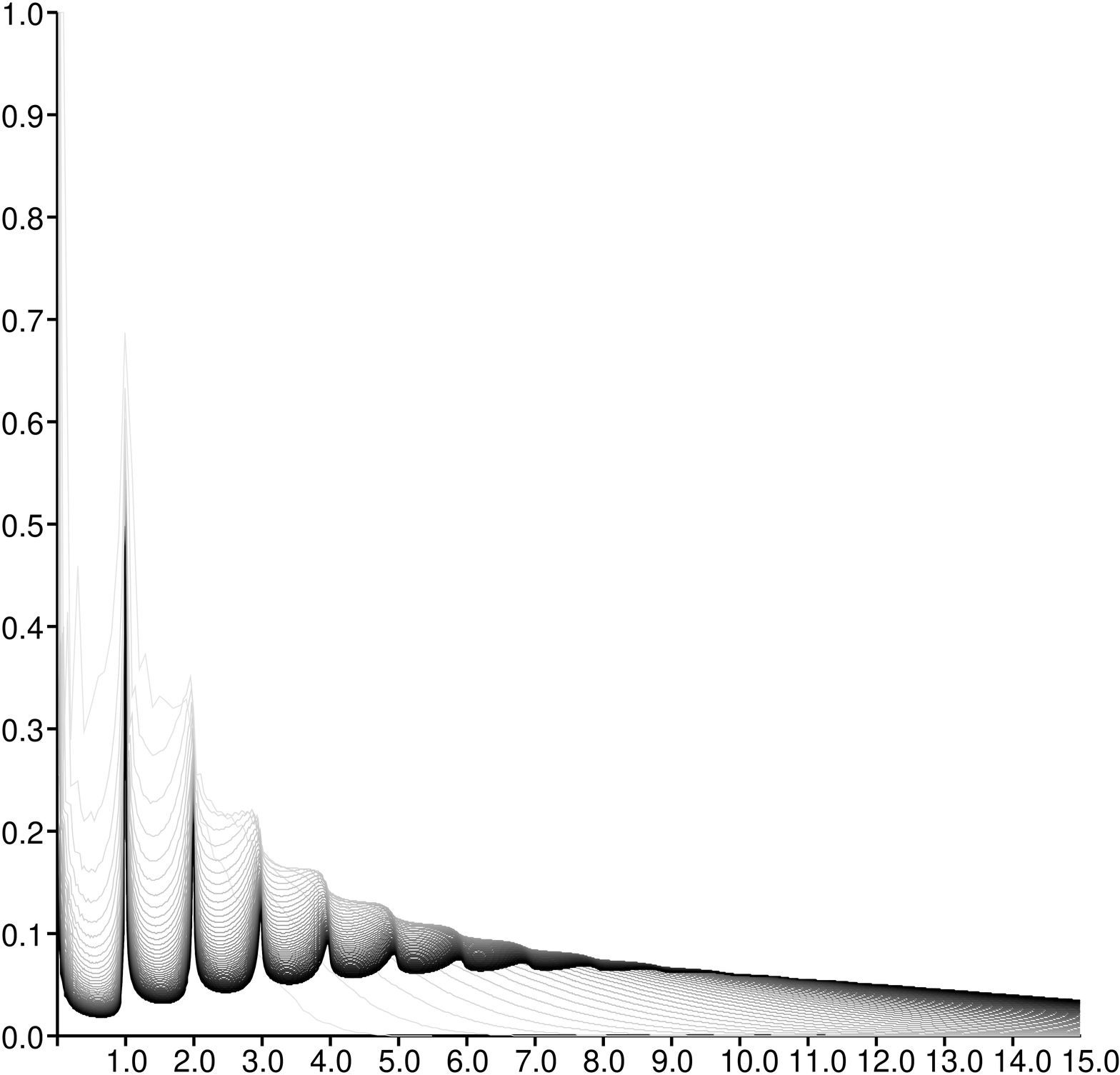}
\includegraphics[width=.45\textwidth]{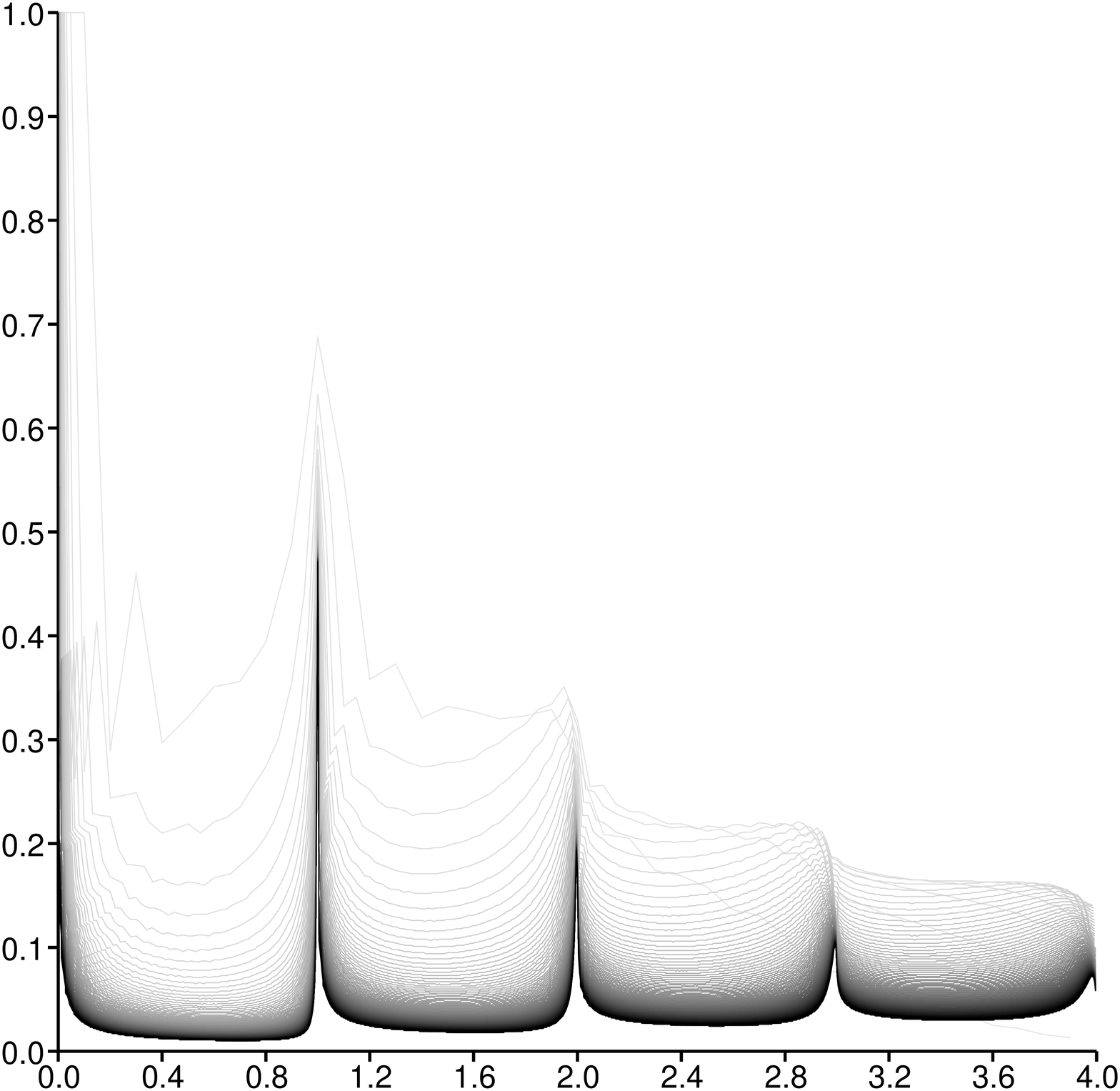}
\caption{Distribution for the first $500$ and $1000$  polynomials}
\label{fig:des1}\label{fig:des2}
\end{center}
\end{figure}

\section{Conclusion}
The Figures (\ref{fig:des1}) and (\ref{fig:des2}) point out some convergence of the renormalized distribution. A question arises naturally:
\begin{itemize}
\item Is it possible to retrieve the asymptotic of the coefficients from their closed form ?
\end{itemize}
Note first that the peaks that appear at integer abscissa $x, x \ll \sqrt{n}$ do not converge to $0$. 
In fact consider the trees made of a path $\pi_0=0,\pi_1,\ldots,\pi_x$ of length $x$ and a plane tree of size $n-x$ rooted in $\pi_x$. The label of $\pi_x$ is $n+(n-1)+(n-2)+\ldots+(n-x+1)$ thus tends to $x$ after the renormalization by $n$. Moreover there are $C_{n-x-1}$ such trees. Since $p_1 = C_{n}$ the peak is higher than $C_{n-x-1}/C_n \sim (1/4)^{x-1}$.

\bibliographystyle{plain}

\bibliography{BibliographieGenerale}

\end{document}